\documentclass[12pt]{article}
\usepackage[utf8]{inputenc} 
\usepackage[english]{babel}
\usepackage{amsmath, amssymb, amsthm, latexsym}
\usepackage{color}
\usepackage{hyperref}

\usepackage{currfile}
\usepackage{navigator}
\IfFileExists{\currfilename}{\embeddedfile{sourcefile}{\currfilename}}{}

\newcommand\highlightForReviewer{\color{black}}

\title{On extended 1-perfect bitrades%
\thanks{
E.A.B. and D.S.K. are with the Sobolev Institute of mathematics,
Novosibirsk 630090, Russia. E-mail: 
bespalovpes@mail.ru, dk@ieee.org
\\
This is the accepted version of the article E.\,A.\,Bespalov, D.\,S.\,Krotov. On extended $1$-perfect bitrades. 
\emph{Discrete Mathematics} 348(1):114138/1--12, 2025, doi
\href{https://doi.org/10.1016/j.disc.2024.114222}{10.1016/j.disc.2024.114222}} 
}
\author{Evgeny A. Bespalov and Denis S. Krotov}
\date{}

\makeatletter
\def\@seccntformat#1{\csname the#1\endcsname.\ } 
\makeatother

\newtheorem{theorem}{Theorem}
\newtheorem{lemma}{Lemma}
\newtheorem{proposition}{Proposition}
\newtheorem{corollary}{Corollary}
\theoremstyle{remark}
\newtheorem{remark}{Remark}
\theoremstyle{definition}
\newtheorem{definition}{Definition}

\newcommand\Cyl[2]{B_{#2}(#1)}
\newcommand\ZZ{{\mathbb Z}}
\newcommand\RR{{\mathbb R}}

\newcommand\prj[2]{{#1^{\lfloor #2 \rfloor}}}
\newcommand\VVH[2]{{\ZZ_{#2}^{#1}}}

\begin{document}

\maketitle

\begin{abstract}
Extended $1$-perfect codes
in the Hamming scheme $H(n,q)$
can be equivalently defined as
codes that turn to $1$-perfect codes 
after puncturing in any coordinate,
as completely regular codes 
with certain intersection array,
as uniformly packed codes 
with certain weight coefficients,
as diameter perfect codes 
with respect to a certain anticode,
as distance-$4$ codes with certain dual distances.
We define extended $1$-perfect bitrades
in $H(n,q)$ in five different manners,
corresponding to the different
definitions
of extended $1$-perfect codes,
and prove the equivalence
of these definitions of extended $1$-perfect bitrades.
 For $q=2^m$, we prove that such bitrades exist if and only if $n=lq+2$. For any $q$,
 we prove the nonexistence of extended $1$-perfect bitrades if $n$ is odd.

\end{abstract}

\def\VV{{\scriptscriptstyle\mathrm{V}}}

\section{Introduction}

In this note, we discuss several different definitions of an extended $1$-perfect bitrade in Hamming graphs $H(n,q)$, $q>2$, 
prove their equivalence,  prove the nonexistence of extended $1$-perfect bitrades for odd $n$ (generalising 
some of the results of~\cite{Bespalov:extmds} on the nonexistence of extended $1$-perfect codes), and prove the existence for $q=2^m$, $n=lq+2$.
As usual, for different combinatorial objects like some classes of codes, designs, orthogonal arrays, etc., 
bitrades are defined as a generalization of the difference pair 
$(C_0\backslash C_1, C_1\backslash C_0)$ of two objects $C_0$ and $C_1$ from the initial class.
The main question is the following: 
which property of the initial class has to be generalized to form the definition of a bitrade.
The five definitions of extended $1$-perfect bitrades considered below reflect different definitions of an extended $1$-perfect
code as 
\begin{itemize}
 \item a code whose projection in any direction 
       (puncturing in any coordinate) 
       gives a $1$-perfect code 
       (this definition reflects 
       a traditional meaning of the concept 
       ``extended code'');
 \item a completely regular code \cite{Neumaier92}, 
         with certain parameters;
 \item a uniformly packed 
       (in the wide sense \cite{BZZ:1974:UPC}) code, 
       with certain parameters;
 \item a diameter-$3$ perfect code \cite{AhlAydKha};
 \item a distance-$4$ code with certain spectral constraints, 
       namely, the characteristic function of the code
       belongs to the sum of three eigenspaces of $H(n,q)$,
 with the eigenvalues $(q-1)n$ 
 (the maximum eigenvalue of $H(n,q)$),
 $q-2$, and $-n$ (the minimum eigenvalue of $H(n,q)$); 
 these spectral restrictions can be also expressed 
 in terms of the components of the dual distance distribution.
\end{itemize}
We note that our definitions are not in agree with known definition of
extended $1$-perfect bitrades in the binary Hamming graph $H(n,2)$.
The reason is that $H(n,2)$, the $n$-cube, is a bipartite graph,
and it is natural to require that the components of the bitrade
lie in one of its parts, either even or odd (which also allows
to consider bitrades as objects in the halved $n$-cube).
Relaxing this restriction leads to some degenerate cases, 
which are not interesting from a theoretical point of view.
Of course, $H(n,q)$ can also be considered as a $q$-partite graph,
but the parts are not uniquely defined, and we cannot guarantee
for an extended $1$-perfect code that it lies in one of them;
neither it is natural to require this from a bitrade.

We finish the introduction by the observation that
extended $1$-perfect codes
do exist in $H(2^k+2,2^k)$ (see, e.g., \cite[\S11.5, Theorem~10]{MWS}), so the theory developed is applicable to the difference pairs of extended $1$-perfect codes, not only to their abstract generalisations.

The next section contains some preliminaries. 
In Section~\ref{s:codes}, we consider
five equivalent definitions of extended $1$-perfect codes in Hamming graphs (four of them are given in the form of propositions).
In Section~\ref{s:def},
we define extended $1$-perfect bitrades
in five different ways.
In Section~\ref{s:equivalence}, 
we show the equivalence of the five definitions of extended $1$-perfect bitrades. In Section~\ref{s:necessary},
we construct extended $1$-perfect bitrades
for any admissible length~$n$ in the case of the alphabet of size~$2^m$ and show
the nonexistence for odd~$n$ in the case of arbitrary alphabet. Section~\ref{s:conc} is concluding.

\section{Preliminaries}

In this paper, graphs are 
unordered, connected, without loops and multiple edges.
A~\emph{code} in a graph 
$G=(V,E)$ is an arbitrary subset 
of~$V$ of cardinality at least~$2$.
The \emph{code distance} of a code~$C$ is the minimum distance between two different vertices of~$C$ 
(the distance $d(x,y)$ in a graph,
by default, is the minimum
length of a path between 
given vertices~$x$ and~$y$).
Given~$n$ and~$q$, 
the \emph{Hamming graph} 
is a graph with the vertex set 
$\mathbb Z^n_q=\{(x_1,\ldots,x_n):x_i \in \mathbb Z_q\}$, 
where two vertices~$x$ and~$y$ are adjacent if and only if they
differ in exactly one position.

The set of vertices
at distance at most~$r$
from a given vertex~$x$
is called the 
\emph{radius-$r$ ball}
with centre~$x$.
The radius-$1$ ball with centre~$x$
will be denoted by~$B(x)$.
The set $S_i(x)=\{y:\ d(x,y)=i\}$
is called the 
\emph{radius-$i$ sphere} 
with centre~$x$.

The next few notations are for
Hamming graphs only, although
the concept of cylinder
can be generalized
to an arbitrary graph and plays
an important role in the generalization
of the Hamming bound on the size
of a code of even code distance.
For $y=(y_1,\ldots,y_n)$ in~$\mathbb Z_q^n$, 
denote by~$y^a_i$, $a\in \mathbb Z_q$, $i\in\{1,\ldots,n+1\}$, the vertex $(y_1,\ldots,y_{i-1},a,y_{i},\ldots,y_n)$ of~$H(n+1,q)$.
Denote $e^a_i=\overline{0}^a_i$, where
$\overline{0}$
is the all-zero word $(0,\ldots,0)$.
Note that $e^0_i$ is the all-zero
word as well,
and for any vertex~$x$
the set $\{x+e^a_i:\ a \in \mathbb Z_q\}$
is the maximum clique in the Hamming graph.
The set $\Cyl{x}{i}=\underset{a \in \mathbb Z_q}\cup B(x+e^a_i)$
is called a \emph{cylinder} (of radius~$1$).

A code $C$ in a graph $G=(V,E)$  is called \emph{$1$-perfect} 
if 
$|B(x) \cap C|=1$ for every~$x$ 
in~$V$.  
The pair $(T_{+},T_{-})$ of disjoint nonempty subsets of~$V$  
is called a \emph{$1$-perfect bitrade} in $H(n,q)$ if 
$|B(x) \cap T_{+}|=|B(x) \cap T_{-}| \le 1$
for every~$x$ in~$V$.  
For two disjoint subsets 
$A_+$ and~$A_-$ of~$V$,
a function $f:V \to \{0,\pm 1\}$ is called the \emph{characteristic function} of the ordered pair $(A_+,A_-)$
if $f(x)=1$ for every $x \in A_+$, $f(x)=-1$ for every $x \in A_-$, and $f(x)=0$ otherwise.

The following fact is well known and straightforward.

\begin{lemma}\label{l:dist3}
Let $G$ be $H(n,q)$, $n>1$, or any other graph
such that $B(z)\not\subseteq B(y)$ for any  different vertices~$y$,~$z$.
If $(T_{+},T_{-})$ is a $1$-perfect bitrade in $G$, then 
$T_{+}$ and $T_{-}$ are codes with code distance~$3$.
\end{lemma}
\begin{proof}
Since $|B(x) \cap T_{+}| \le 1$ by the definition,
the code distance of~$T_{+}$ is \emph{at least~$3$}.

Consider a ball $B(z)$ centred in a vertex~$z$ in~$T_{-}$. By the definition of a $1$-perfect bitrade, for each~$x$ in~$B(z)$, there is $y_x\in T_{-}$ such that $d(x,y_x)\le 1$. 
By the hypothesis on the graph~$G$, $y_x$ cannot be the same for all~$x$ in~$B(z)$. Hence, there are 
adjacent $x'$ and~$x''$ such that $y_{x'}\ne y_{x''}$. Then we have $d(y_{x'},y_{x''})\le 3$,
and the code distance of~$T_{+}$ is \emph{at most~$3$}.
\end{proof}

Given a graph $G=(V,E)$, a function 
on the vertex set~$V$ can be naturally treated 
as a vector from $|V|$ values, indexed by 
elements of~$V$.
A function $f:V \to \RR$
is called an \emph{eigenfunction} of~$G$
with eigenvalue~$\lambda$ if $f \not\equiv 0$ and
for every vertex~$x$
\begin{equation}\label{eq:eig}
 \lambda f(x)=\sum_{y:\,\{x,y\}\in E}f(y)
\end{equation}
(i.e., $f$ corresponds to an eigenvector
of the adjacency matrix of~$G$).
The Hamming graph $H(n,q)$ has the following eigenvalues (see, e.g., \cite[Theorem~9.2.1]{Brouwer}): 
$$
\lambda_i(n,q)=n(q-1)-iq, \quad i=0,1,\ldots,n;$$
we will refer to $\lambda_i(n,q)$
as the $i$th eigenvalue of~$H(n,q)$
and denote by $U_i(n,q)$ 
the corresponding eigenspace.

{\highlightForReviewer
For a code $C$ in $H(n,q)$
and $i\in\{1,\ldots,n\}$,
the \emph{$i$-projection} of~$C$ (or simply, a \emph{projection})
is the code in $H(n-1,q)$ obtained
by removing the $i$th entry 
in all elements of~$C$.
For the pair $(T_{+},T_{-})$ of disjoint sets of vertices of~$H(n,q)$,
we define 
the \emph{\mbox{$i$-}projection}
as the pair
$(A_{+} \backslash A_{-},
  A_{-} \backslash A_{+})$,
where 
$A_{+}$ and~$A_{-}$ are the $i$-projections of~$T_{+}$ and~$T_{-}$ respectively.}

For a function
$f:\ \mathbb Z_q^n \to \RR$,
we define its $i$-projection $\prj{f}{i}:\ \mathbb Z_q^{n-1} \to \RR$ in the following way:
$\prj{f}{i} (x)=\underset{y=x^a_i:\  a \in \mathbb Z_q}\sum f(y)$.
Note that the characteristic function 
of the $i$-projection 
of a code~$C$ with distance at least~$2$ 
(or a bitrade $(T_{+},T_{-})$) is 
$\prj{f}{i}$, where $f$ is the characteristic function of the code (respectively, bitrade).  

\begin{lemma}\label{l:fiproekt}
\highlightForReviewer
Let $f$ be an eigenfunction of $H(n,q)$ corresponding to the $j$th eigenvalue.
\begin{enumerate}
\item \cite{ValVor:Hamm} If $j=n$, then 
all projections of~$f$ are (constantly) zero.
\item \cite{ValVor:Hamm} If $j < n$, then all non-zero projections of~$f$ are eigenfunctions of~$H(n-1,q)$ with the $j$th eigenvalue. 
\item If $j < n$, 
then there is a non-zero projection of~$f$. 
\end{enumerate}
\end{lemma}
\begin{proof}[Proof of p.\,3]
Assume $f$ is an eigenfunction. 
Denote by $c_{x,i}$ the set 
$\{x+e^a_i:a \in \mathbb Z_q\}$ 
of vertices. 
If $\prj{f}{i} \equiv 0$ for any $i$, 
then for any~$x$ and~$i$
we have $\sum_{y\in c_{x,i}} f(y)=0$
and
$$
\sum_{y:\,d(x,y)=1} f(y)
= \sum_{i=1}^n
\Big(\sum_{y\in c_{x,i}} f(y) - f(x) \Big) = - n f(x).
$$
The last implies $f\in U_n(n,q)$,
 contradicting the hypothesis of~p.3 of the lemma.
\end{proof}

\section{Extended 1-perfect codes}\label{s:codes}

Here we consider five equivalent definitions
of extended $1$-perfect codes. 
These definitions are not used 
in the rest of the paper, but they 
motivate the corresponding 
five definitions of extended $1$-perfect
bitrades.

\subsection{Projections}\label{s:ext-proj}
The extension is usually considered as the operation 
inverse to the projection
(also known as \emph{puncturing} 
in coding theory).
This leads to the following
definition of an
extended $1$-perfect code.

\setcounter{definition}{-1}
\begin{definition}
A code $C$ in $H(n,q)$ is an 
\emph{extended $1$-perfect code}
if every projection of~$C$ is a $1$-perfect code.
\end{definition}
This definition will be considered as the main one,
while the others  
will be given in the form
of propositions.

\subsection{Diameter perfect codes}\label{s:ext-diam}
A distance-$d$ code $C$ in $H(n,q)$ (or any other transitive or distance-regular graph) is called \emph{diameter perfect}
if $|C|\cdot|D|=|\VVH{n}{q}|$ for some set of vertices (anticode)
of diameter $d-1$. The next definition follows from the fact 
that an extended $1$-perfect code is diameter perfect, 
see~\cite{AhlAydKha}, with a cylinder $\Cyl{x}{i}$ as an anticode.

\begin{proposition}
A code $C$ in $H(n,q)$ 
is an \emph{extended $1$-perfect code}
if and only if $|C \cap \Cyl{x}{i}|=1$ 
for every vertex $x$ 
and coordinate $i$.
\end{proposition}
\begin{proof}
Let $C$ be an extended 
$1$-perfect code in~$H(n,q)$.
The cardinality of~$C$ equals the cardinality of a perfect code in $H(n-1,q)$,
i.e., $q^{n-1}/(1+(q-1)(n-1))$.
The cardinality 
of any cylinder~$D$
in $H(n,q)$ equals $q\cdot(1+(q-1)(n-1))$.
We see that
\begin{equation}\label{eq:cyl-anti}
   |C|\cdot |D| = q^n=|\VVH{n}{q}|.  
\end{equation}
The code distance of~$C$ is~$4$,
because the code distance 
of each its projection is~$3$.
Since a cylinder has diameter~$3$,
it cannot contain more than~$1$
element  of~$C$. 
Denote by $N$ the total number of cylinders. 
Clearly, each vertex of $H(n,q)$
is contained in the same number $N\cdot |D| / q^n$ of cylinders. Hence, there are
$|C|\cdot N\cdot |D| / q^n$
cylinders intersecting with~$C$.
Utilising \eqref{eq:cyl-anti},
we see that this number is~$N$, 
i.e., each cylinder intersects with~$C$.

Inversely, if each cylinder intersects 
with~$C$ in exactly one element,
then for every~$i$ the $i$-projection
of~$C$ intersects with each radius-$1$
ball in exactly one element
(indeed, every such ball $B(x)$ in $H(n-1,q)$
is the $i$-projection
of the cylinder $B(x_i^0)$).
In this case, by the definitions,
all projections of~$C$ are 
$1$-perfect codes and hence~$C$ is 
an extended $1$-perfect code.
\end{proof}

\subsection{Completely regular codes}\label{s:ext-cr}
For a code~$C$ in a connected graph $G=(V,E)$, 
denote by~$C^{(i)}$ the set of vertices
at distance~$i$ from~$C$.
The \emph{covering radius} of~$C$, $\rho(C)$,
is the maximum~$i$ such that $|C^{(i)}|>0$. 
A~code~$C$ is called 
\emph{completely regular} if there are constants
$b_i$, $a_i$, $c_{i}$,
$i=0,1,\ldots,\rho(C)$
(\emph{intersection numbers})
such that every vertex in~$C^{(i)}$ has exactly~$b_i$ neighbours in $C^{(i+1)}$, 
$a_i$ neighbours in~$C^{(i)}$, 
and~$c_{i}$ neighbours in $C^{(i-1)}$, 
$i=0,1,\ldots,\rho(C)$.
Straightforwardly, 
\begin{equation}\label{eq:AFFS}
    AF=FS,
\end{equation}
where~$A$ is the adjacency matrix of the graph,
$F$ is the $|V|\times (\rho(C)+1)$ matrix 
whose columns are the characteristic functions of 
$C$, $C^{(1)}$, \ldots, $C^{(\rho(C))}$, respectively, and
$S$ is the tridiagonal matrix with 
main diagonal
$(a_0$, $a_1$, \ldots, $a_{\rho(C)})$,
superdiagonal
$(b_0$, $b_1$, \ldots, $b_{\rho(C)-1})$,
and subdiagonal 
$(c_1$, $c_2$, \ldots, $c_{\rho(C)})$.
Note that if the graph is regular of degree~$d$, then 
$a_i=d-b_i-c_i$ for every~$i$.

\begin{proposition}[{see, e.g., \cite[(F.2)]{BRZ:CR}, \cite{Bespalov:extmds}}]
A code $C$ in $H(n,q)$ 
is an \emph{extended $1$-perfect code}
if and only if it is completely regular with covering radius~$2$ and
intersection numbers $b_0=n(q-1)$, $b_1=(n-1)(q-1)$, $c_1=1$, $c_2=n$.
\end{proposition}

\subsection{Uniformly packed codes}\label{s:ext-up}
A code $C$ in $H(n,q)$ is called \emph{uniformly packed}
in the wide sense \cite{BZZ:1974:UPC}
if there are numbers 
$\alpha_0$, \ldots, $\alpha_{\rho(C)}$ such that
for every $x$ from $\VVH{n}{q}$ it holds
$$\alpha_0f_0(x)+ \ldots  +\alpha_{\rho(C)}f_{\rho(C)}(x)=1$$
where $f_i(x)$ is the number 
of elements of~$C$ at distance~$i$ 
from~$x$.

\begin{remark}
The coefficients $\alpha_i$ are not required to be all non-negative.
For example, \cite[(14)]{BZZ:1974:UPC}
the direct product of two copies of the same $1$-perfect code in $H(\frac n2,q)$ is uniformly packed
in $H(n,q)$
with $\alpha_0=\alpha_2=1$ and $\alpha_1 = 1-\frac14(q-1)(n-2)$.
\end{remark}

\begin{proposition} A code $C$ in $H(n,q)$ is an \emph{extended $1$-perfect code}
if and only if it is uniformly packed in the wide sense with 
$\rho(C)=2$, 
$\alpha_0 = \alpha_1=1$,
$\alpha_2=\frac2n$.
\end{proposition}

\subsection{Codes with two dual distances}
\label{s:ext-dual}
As follows from linear algebra, the characteristic 
function $\chi_C$ of a code $C$ in $H(n,q)$ 
is uniquely represented as the sum
$$ \chi_C = \phi_0+\ldots+\phi_n$$
where $\phi_i$ is an eigenfunction with eigenvalue $\lambda_i(n,q)$
or the constant zero.
The indices~$i$ such that $\phi_i\not\equiv 0$
are called the \emph{dual distances} of $C$
(essentially, they are the distances in the dual
distance distribution obtained by the MacWilliams transform,
see e.g.~\cite{MWS}).

\begin{proposition}  A code $C$ in $H(n,q)$ 
is an extended $1$-perfect code 
if it has code distance at least $4$ and dual distances 
$0$, $n-\frac{n-2}q -1$, and $n$. 
\end{proposition}

\begin{proof}
Assume that $C$ is a distance-$4$ code with
dual distances
$0$, $n-\frac{n-2}q -1$, and~$n$. 
For every $i$, the $i$-projection of~$C$ 
has the code distance 
at least~$3$ and by Lemma~\ref{l:fiproekt}(1,2)
its characteristic function
has eigenvalues~$-1$ and~$n(q-1)$. Hence, 
this projection is a $1$-perfect code,
and $C$ is an extended $1$-perfect code
by the definition.

On the other hand, 
an extended $1$-perfect code
has eigenvalues 
$n(q-1)$, $q-2$, and $-n$,
see e.g.~\cite[(F.2)]{BRZ:CR}, ~\cite{Bespalov:extmds}.
\end{proof}

\section{Definitions of extended 1-perfect bitrades}\label{s:def}
We are now ready to introduce
five definitions of extended $1$-perfect bitrades, 
which correspond
to the five definitions
of extended $1$-perfect codes
from Sections~\ref{s:ext-proj}--\ref{s:ext-dual}, respectively.
\begin{definition}[using projection]\label{d:proection}\label{d:first}
A pair $(T_{+},T_{-})$ of disjoint subsets 
of~$\VVH{n+1}{q}$
is called an \emph{extended $1$-perfect bitrade} 
in~$H(n+1,q)$ 
if $T_{+}$ and~$T_{-}$ are codes with code distance~$4$ and 
 for every $i \in \{1,\ldots,n\}$ the $i$-projection of $(T_{+},T_{-})$
is a $1$-perfect bitrade in~$H(n,q)$.
\end{definition}

\begin{definition}[diameter-perfect bitrade, using anticodes]\label{d:bochka}
Given a collection~$\mathcal A$
of diameter-$d$
sets (anticodes)
of vertices of a graph,
a pair $(T_{+},T_{-})$ of disjoint 
sets of vertices is called 
a \emph{diameter-$d$-perfect bitrade}
(with respect to~$\mathcal A$)
if for every $A\in\mathcal A$
it holds 
$|A \cap T_{+}|=|A \cap T_{-}| \le 1$.
A pair $(T_{+},T_{-})$ of disjoint subsets of~$\VVH{n}{q}$ is called an \emph{extended $1$-perfect bitrade} in $H(n,q)$
if it is a diameter-$d$-perfect bitrade
with respect to the set of cylinders
$\Cyl{x}{i}$, $x \in \VVH{n}{q}$, $i \in \{1,\ldots,n\}$.
\end{definition}

\begin{lemma}\label{l:dist44}
Let $(T_{+},T_{-})$ be an extended $1$-perfect bitrade according to Definition~\ref{d:bochka}.
Then $T_{+}$ and $T_{-}$ are codes with code distance~$4$.
\end{lemma}
\begin{proof}
Seeking a contradiction, we assume that there are two vertices~$x$ and~$y$ in~$T_{+}$ such that $d(x,y) \le 3$. There is a path $(x, a, b, y)$ 
or $(x=a, b, y)$. Then $x,y \in \Cyl{a}{i}$, where $i$ is the position in which $b$ and $y$ differ, and we get 
$|\Cyl{a}{i} \cap T_{+}| \ge 2$, contradicting 
$|\Cyl{a}{i} \cap T_{+}| \le 1$ 
from Definition~\ref{d:bochka}
of bitrade. For $T_{-}$, the proof is similar.     
\end{proof}

\begin{definition}[completely regular bitrade, using matrix equation]\label{d:matrixeq}
A pair $(T_{+},T_{-})$ of disjoint sets of vertices of a graph of order~$N$
is called an 
\emph{completely regular bitrade} 
with a $k\times k$ tridiagonal quotient matrix~$S$
if there is an $N \times k$ matrix~$F$  
whose columns are indexed 
by the vertices of the graph
such that:
\begin{itemize}
    \item[(i)] the first column of~$F$ is the characteristic function of $(T_{+},T_{-})$;
    \item[(ii)] 
    in each row is either all-zero or it contains exactly two non-zero elements, one~$1$ and one~$-1$;
    \item[(iii)] it holds $AF=FS$, where $A$ is the adjacency matrix of the graph. 
\end{itemize} 
A pair $(T_{+},T_{-})$ of disjoint subsets of $\VVH{n}{q}$ is called an \emph{extended $1$-perfect bitrade}
in $H(n,q)$ if
$T_{+}$ and $T_{-}$ are codes 
with code distance~$4$
and
$(T_{+},T_{-})$ is a 
completely regular bitrade 
with quotient matrix
$$
S=\left(\begin{array}{c@{\ \ \ }cc}
0 & n(q-1) & 0 \\
1 & q-2 & (n-1)(q-1) \\ 
0 & n & n(q-2)
\end{array}\right).$$
\end{definition}

{\highlightForReviewer
As follows from equation 
\eqref{eq:AFFS} for completely regular codes, the difference pair
$(C_1\backslash C_2,C_2\backslash C_1)$ of two different completely regular codes with the same quotient matrix is a completely regular bitrade.
For a completely regular code, the incidence matrix~$F$ is defined
in a unique way.
For a completely regular bitrade,
the uniqueness of~$F$ is not 
evident from the definition,
but can be proved separately.
}

\begin{lemma}
If $(T_{+},T_{-})$ is a completely
regular bitrade with quotient matrix~$S$, where~$S$ is a tridiagonal matrix with no 
zeros in the subdiagonal,
then the matrix~$F$
satisfying {\rm(i)--(iii)}
is unique.
\end{lemma}
\begin{proof}
Assume that two matrices~$F_1$ and~$F_2$ 
satisfy the definition; denote their difference by~$F'$ 
and its columns 
by~$x_0$, $x_1$, \ldots, $x_{k-1}$.
From~(iii) we have $AF'=F'S$;
from~(i), 
the first column~$x_0$ 
contains only zeros 
(because the first columns of~$F_1$ and~$F_2$ 
are characteristic functions 
of the same bitrade).
Denoting the elements of~$S$
by $S_{i,j}$, $i,j=0,\ldots,k-1$,
we see that the first 
column of the right part
in $AF'=F'S$ is~$S_{0,0}x_0+S_{1,0}x_1$,
while the first 
column of the left part
is~$Ax_0$. Since $x_0$ is all-zero
and $S_{1,0}\ne 0$,
we conclude that~$x_1$ 
is all-zero too.
The same is true for~$x_2$,
because from~(iii)
we have 
$Ax_1=S_{0,1}x_0+S_{1,1}x_1+S_{2,1}x_2$, $S_{2,1}=0$, then for~$x_3$
(if any), and so on.
We conclude that~$F'$ is all-zero and
$F_1=F_2$.
\end{proof}

For an arbitrary set~$T$ of vertices of a graph,
a vertex~$x$, and a collection
$\bar\alpha = (\alpha_0,\ldots,\alpha_\rho)$,
define the \emph{$\bar\alpha$-weight}
$w_{T,\bar\alpha}(x)$ of~$x$
with respect to~$T$ in the following way:
\begin{equation}
    \label{eq:wT} \highlightForReviewer
w_{T,\bar\alpha}(x)=\sum_{j=0}^\rho\alpha_j\cdot|S_j(x) \cap T|.
\end{equation}

\begin{definition}[uniformly packed bitrade, using weighted balls]\label{d:weighball}
A~pair $(T_{+},T_{-})$ of disjoint sets of 
vertices of a graph is called a 
\emph{uniformly packed bitrade}
with coefficients $\bar\alpha = (\alpha_0,\ldots,\alpha_\rho)$, $\alpha_\rho>0$,
if for any vertex~$x$ it holds $w_{T_{+},\bar\alpha}(x)=w_{T_{-},\bar\alpha}(x) \le 1$.
A uniformly packed bitrade in $H(n,q)$ with coefficients $\bar\alpha = (1,1,\frac2n)$ is called an \emph{extended $1$-perfect bitrade}. 
\end{definition}

For subsets $T_{+}$ and $T_{-}$ of $\VVH{n}{q}$,
we slightly simplify the notation for the 
$(1,1,\frac2n)$-weight with respect to~$T_+$ and~$T_-$:
\begin{equation}
    \label{eq:w+}\highlightForReviewer
w_{+}(x)=w_{T_{+},\bar\alpha}(x)
, \quad w_{-}(x)=w_{T_{-},\bar\alpha}(x), \qquad \bar\alpha=(1,1,\textstyle\frac2n).
\end{equation}

\begin{lemma}\label{l:dist4}
Let $(T_{+},T_{-})$ be an extended $1$-perfect bitrade according to Definition~\ref{d:weighball}.
Then $T_{+}$ and $T_{-}$ are codes with code distance~$4$.
\end{lemma}
\begin{proof}
Seeking a contradiction,
assume that 
$1\le d(z,y) \le 3$
for some
$z$ and~$y$ from~$T_{+}$ (for $T_{-}$, the proof is similar).
In $H(n,q)$, there is a 
vertex~$x$ at distance~$1$ from~$z$ and 
at distance
$d(x,y)-1$ from~$y$;
{\highlightForReviewer
that is, $z\in S_1(x) \cap T_+$ and $y\in S_{d(x,y)-1}(x) \cap T_+$.
From the definition of $w_{+}$ in~\eqref{eq:w+} and~\eqref{eq:wT}, 
we see that }
\begin{multline*}
 \label{eq:wT} \highlightForReviewer
    w_{+}(x) =
w_{T_+,\bar\alpha}(x)=\sum_{j=0}^2\alpha_j\cdot|S_j(x) \cap T_+|
\\ \highlightForReviewer
\ge
\alpha_1\cdot |\{z\}| +
\alpha_{d(x,y)-1}\cdot |\{y\}|
= \alpha_1 + \alpha_{d(x,y)-1}
\ge 1 + \frac{2}{n}>1,
\end{multline*}
which contradicts 
{\highlightForReviewer the condition $w_{T_+,\bar\alpha}(x)\le 1$}
in Definition~\ref{d:weighball}.
\end{proof}

\begin{definition}[using eigenspaces]\label{d:eigensubspace}\label{d:last}
A pair $(T_{+},T_{-})$ of disjoint subsets of $\VVH{n}{q}$ is called an \emph{extended $1$-perfect bitrade} in $H(n,q)$
if $T_{+}$ and $T_{-}$ are codes with code distance~$4$ and 
the characteristic function~$f$ of the pair $(T_{+},T_{-})$
belongs to the direct sum
of eigenspaces corresponding to the eigenvalues $-n$ and $q-2$,
that is,
$$\highlightForReviewer
f \in U_n(n,q) \oplus U_l(n,q),\quad
\mbox{where }
l=\frac{n(q-1)-q+2}{q}.
$$
\end{definition}

\section{Equivalence of the definitions}\label{s:equivalence}

One of the main results
of our paper is
the following theorem.

\begin{theorem}
\label{th:equiv}
Definitions~\ref{d:first}--\ref{d:last} are equivalent.
\end{theorem}

We divide the proof into
three parts, Propositions~\ref{p:235}--\ref{p:24}, 
preceding them by a technical lemma about the function~$w_+$, defined in~\eqref{eq:wT}--\eqref{eq:w+}.
Note that all definitions of bitrades are symmetric with respect to~$T_{+}$ and~$T_{-}$, i.e.,
$(T_{+},T_{-})$ is a bitrade
if and only if $(T_{-},T_{+})$ is.
We will often use this fact
in the proofs to consider only one
of the two cases if the proof itself
is not symmetric with respect to~$T_{+}$ and~$T_{-}$.

For vertices $x$ and $y$ of $H(n,q)$, denote by $\Delta(x,y)$ the set of positions for which $x$ and $y$ differ. 
The cardinality of $\Delta(x,y)$ equals $d(x,y)$.

\begin{lemma}\label{l:weight1}
Let $T_+$ be a code in $H(n,q)$ with code distance~$4$. For any vertex~$x$ of~$H(n,q)$,
the following assertions hold:
\begin{enumerate}
    \item $w_+(x) \le 1;$
    \item if $w_+(x)=1$, then
    \begin{itemize}
        \item 
    either $|B(x) \cap T_+|=1$ and $|S_2(x) \cap T_+|=0$,
    \item
    or $|B(x) \cap T_+|=0$, $|S_2(x) \cap T_+|=n/2$, and $\underset{v \in S_2(x) \cap T_+}\cup \Delta(x,v)=\{1,\ldots,n\}$,
    $\Delta(x,y) \cap \Delta(x,z)=\emptyset$ for different~$y$ and~$z$ from~$S_2(x) \cap T_+$.
    \end{itemize}
\end{enumerate}
\end{lemma}
\begin{proof}
If there is a vertex in $B(x) \cap T_+$, 
then by the code distance $|B(x) \cap T_+|=1$, $|S_2(x) \cap T_+|=0$, and $w_+(x)=1$.

Now assume $B(x) \cap T_+ = \emptyset$.
Let there are $k$ vertices in $S_2(x) \cap T_+$.
Since the code distance is~$4$, we have $\Delta(x,y) \cap \Delta(x,z) = \emptyset$ for any different~$y$ and~$z$ in $S_2(x) \cap T_+$.
It follows that $2k \le n$, 
$k \le n/2$, and 
$w_+(x) \le 1$.
Moreover, if $w_+(x) = 1$, 
then $k=n/2$ and 
$|\underset{v \in S_2(x) \cap T_+}\cup \Delta(x,v)| = \frac{n}{2}\cdot 2 = |\{1,\ldots,n\}|$.
\end{proof}

\begin{proposition}\label{p:235}
Definitions~\ref{d:proection},~\ref{d:weighball}, and \ref{d:bochka} are equivalent.
\end{proposition}
\begin{proof}
$\ref{d:bochka} \Longrightarrow \ref{d:proection}$.
Let $(T_{+},T_{-})$ be an extended $1$-perfect bitrade in $H(n,q)$ according to Definition~\ref{d:bochka}, and
let $i \in \{1,\ldots,n\}$ be an arbitrary position. 
Consider the $i$-projection of the bitrade, denote it by $(T'_{+},T'_{-})$.
Let $x$ be an arbitrary vertex of $H(n-1,q)$ and let $y=x^0_i$.
If there are no vertices from $(T_{+} \cup T_{-})$ in $\Cyl{y}{i}$, 
then there are no vertices from $(T'_{+} \cup T'_{-})$ in $B(x)$.
Otherwise, let $\Cyl{y}{i} \cap T_{+} =\{v\}$, $\Cyl{y}{i} \cap T_{-}=\{u\}$.
If $u$ and $v$ differ only in position~$i$, 
then there are no vertices from $(T'_{+} \cup T'_{-})$ in~$B(x)$; otherwise, 
we have $|B(x) \cap T'_{+}|=|B(x) \cap T'_{-}|=1$.

$\ref{d:weighball} \Longrightarrow \ref{d:bochka}$.
Let $(T_{+},T_{-})$ be an extended $1$-perfect bitrade in $H(n,q)$ according to Definition~\ref{d:weighball}.
Consider $\Cyl{x}{i}$ for an arbitrary $x \in \VVH{n}{q}$ and $i \in \{1,\ldots,n\}$.
Since $T_{+}$ and $T_{-}$ are codes with code distance $4$ (by Lemma~\ref{l:dist4}), there 
are at most one vertex from~$T_{+}$ and at most one vertex from~$T_{-}$ in $\Cyl{x}{i}$.
It is sufficient to prove that if there is a vertex~$y$ from~$T_{+}$  in $\Cyl{x}{i}$, 
then there is a vertex from~$T_{-}$ in $\Cyl{x}{i}$.

%
%
Since $y \in \Cyl{x}{i}$, there is a vertex~$z$
at distance~$1$ from~$y$
such that $z=x+e^a_i$ for some $a \in \mathbb Z_q$; since $y \in T_{+}$, we have $w_{+}(z)=1$.
By Lemma~\ref{l:weight1}, we have that either $z$ is adjacent to some vertex~$u$ from~$T_{-}$, 
or there are $n/2$ vertices from~$T_{-}$ at the distance~$2$ from~$z$.
In the first case, $u$ belongs to~$\Cyl{x}{i}$.
In the second case, by Lemma~\ref{l:weight1} there is a vertex $v \in T_{-}$ at the distance~$2$ from~$z$ such that 
$z$ and $v$ differ in position~$i$ and some other position; hence, $v$ belongs to $\Cyl{x}{i}$.

$\ref{d:proection} \Longrightarrow \ref{d:weighball}$.
Let $(T_{+},T_{-})$ be an extended $1$-perfect bitrade in~$H(n,q)$ according to Definition~\ref{d:proection}.
Consider an arbitrary vertex~$x$. 
Let us prove that $w_{+}(x)=w_{-}(x)$.
We consider four cases. 

(a) Suppose, $x \in T_{+}$ (similarly, $x \in T_{-}$). 
In this case, $w_{+}(x)=1$. 
If there is a vertex from~$T_{-}$ in the neighbourhood of~$x$, then $w_{-}(x)=1$.
Otherwise, let $y_1,\ldots,y_k$ be the elements of~$T_{-}$ at distance~$2$ from~$x$.
Since the code distance is~$4$, the sets $\Delta(x,y_1),\ldots,\Delta(x,y_k)$ are pairwise disjoint. 
If there is~$i$ in $\{1,\ldots,n\}$
that does not belong to the union of these sets,
then 
the vertex that corresponds to~$x$ in the $i$-projection has no neighbours in the projection of~$T_-$; 
this contradicts the definition of $1$-perfect bitrade.
Therefore, $k=n/2$ and $w_{-}(x)=1$.

(b) Suppose $x \not\in (T_{+} \cup T_{-})$ and $|S_1(x) \cap T_{+}|=|S_1(x) \cap T_{-}|=1$.
Obviously, $w_{+}(x)=w_{-}(y)=1$ in this case.

(c) Suppose $x \not\in (T_{+} \cup T_{-})$, $|S_1(x) \cap T_{+}|=1$, and $|S_1(x) \cap T_{-}|=0$ (similarly, $|S_1(x) \cap T_{+}|=0$, $|S_1(x) \cap T_{-}|=1$). 
Analogously to (a), there are $n/2$ vertices in~$T_{-}$ at distance~$2$ from~$x$; hence, $w_{-}(x)=1$.

(d) Suppose $x \not\in (T_{+} \cup T_{-})$ and $|S_1(x) \cap T_{+}|=|S_1(x) \cap T_{-}|=0$.
Using arguments similar to~(a) and~(c), we will prove that in~$S_2(x)$,
the number of vertices from~$T_{+}$ 
equals 
the number of vertices from~$T_{-}$.
Seeking a contradiction, assume
that there is~$i$ in
$\underset{v \in S_2(x) \cap T_{+}}\cup \Delta(x,v) 
\setminus
\underset{v \in S_2(x) \cap T_{-}}\cup \Delta(x,v)$.
Consider the $i$-projection of the bitrade, the pair $(T'_{+},T'_{-})$.
We see that $B(x')$, where $x'$ is the $i$-projection image
of~$x$, 
contains a vertex from~$T'_{+}$ but
does not contain a vertex from~$T'_{-}$, 
a contradiction.
Therefore, $|S_2(x) \cap T_{+}|=|S_2(x) \cap T_{-}|$ and $w_{+}(x)=w_{-}(x)$.

In all the cases, we get $w_{+}(x)=w_{-}(x)$; hence,
$(T_{+},T_{-})$ is an extended $1$-perfect bitrade according to Definition~\ref{d:weighball}.
\end{proof}

\begin{proposition}\label{13}
 Definitions~\ref{d:matrixeq} and \ref{d:weighball} are equivalent.
\end{proposition}
\begin{proof}
$\ref{d:weighball} \Longrightarrow \ref{d:matrixeq}$.
Let $(T_{+},T_{-})$ be an extended $1$-perfect bitrade in $H(n,q)$ according to Definition~\ref{d:weighball}. 
Let us construct a matrix~$F$ that satisfies Definition~\ref{d:matrixeq}. 
The first column of $F$ is the characteristic function of $(T_{+},T_{-})$.
For a vertex~$v$, denote by $f(v)$ the $v$th row of~$F$. We define the values of $f(\cdot)$ as follows:
\begin{itemize}
    \item for $x \in T_{+}$,
    $f(x)=
    \left\{
    \begin{array}{l@{\,}r@{\,}r@{\,~}r@{\,}ll}
    (&+1,&-1,&0&)& \mbox{if $S_1(x)\cap T_{-} \ne \emptyset$,} \\
    (&+1,&0,&-1&)& \mbox{if $S_1(x)\cap T_{-} = \emptyset$,}
    \end{array}\right.
    $
    \item for $x \in T_{-}$,
    $f(x)=
    \left\{
    \begin{array}{l@{\,}r@{\,}r@{\,}r@{\,}ll}
    (&-1,&+1,&0&)& \mbox{if $S_1(x)\cap T_{+} \ne \emptyset$,} \\
    (&-1,&0,&+1&)& \mbox{if $S_1(x)\cap T_{+} = \emptyset$,}
    \end{array}\right.
    $
    \item for $x \not\in T_{+}\cup T_{-}$,
    $f(x)=
    \left\{
    \begin{array}{l@{\,}r@{\,}r@{\,}r@{\,}ll}
    (&0,&0,&0&)& \mbox{if $S_1(x)\cap T_{+} \ne \emptyset \ne S_1(x)\cap T_{-}$,} \\
    (&0,&+1,&-1&)& \mbox{if $S_1(x)\cap T_{+} \ne \emptyset=S_1(x)\cap T_{-}$,} \\
    (&0,&-1,&+1&)& \mbox{if $S_1(x)\cap T_{+} = \emptyset \ne S_1(x)\cap T_{-}$,} \\
    (&0,&0,&0&)& \mbox{if $S_1(x)\cap T_{+} = \emptyset = S_1(x)\cap T_{-}$.}
    \end{array}\right.
    $
\end{itemize}
We have defined the matrix~$F$, and
it remains to prove that $AF=FS$.
Let $x$ be an arbitrary vertex of $H(n,q)$.
The neighbourhood of~$x$ is the union of $n$ 
pairwise disjoint cliques of cardinality $q-1$; denote them by $N_1,\ldots,N_n$, i.e., $S_1(x)=\cup_{i=1}^n N_i$.
Consider four cases depending on~$f(x)$.

(a) Let $f(x)=(+1,-1,0)$. 
By definition, $x \in T_{+}$,
there is a vertex $y \in T_{-}$ in~$S_1(x)$, and 
there are no vertices from $(T_{+} \cup T_{-})$ in $S_2(x)$.
The vertex $y$ belongs to $N_i$ for some~$i$,
and $f(y)=(-1,+1,0)$. 
In the remaining $q-2$ vertices of $N_i$, $f$ has value $(0,0,0)$, 
and in the remaining $(n-1)(q-1)$ vertices 
of the neighbourhood, value $(0,+1,-1)$.
Hence, $$AF_x=(-1,+1,0)+(n-1)(q-1) \cdot (0,+1,-1)=\big(-1,\ (n-1)(q-1)+1,\ -(n-1)(q-1)\big),$$ where 
$AF_x$ is the $x$-row of the matrix~$AF$.
On the other hand, $(+1,-1,0)\cdot S$ 
also equals $(-1,(n-1)(q-1)+1,-(n-1)(q-1))$.

(b) Let $f(x)=(+1,0,-1)$. 
By definition, $x \in T_{+}$,
there are no vertices from~$T_{-}$ in~$S_1(x)$, and 
there are $n/2$ vertices from~$T_{-}$ in~$S_2(x)$.
It follows that in each clique~$N_i$, there is exactly one vertex adjacent to some vertex from~$T_{-}$.
Therefore, in the neighbourhood of~$x$,
there are $n$ vertices with value $(0,0,0)$ of~$f$
and $n(q-2)$ vertices with value $(0,+1,-1)$.
We find $AF_x=n(q-2) \cdot (0,+1,-1)$. On the other hand, $(+1,0,-1)\cdot S = (0,n(q-2),-n(q-2))$.

(c) Let $f(x)=(0,+1,-1)$. 
We have $(0,+1,-1)\cdot S=(+1, -n+q-2, n-q+1)$.
By definition, $x \not\in (T_{+} \cup T_{-})$,
there is a vertex $z \in T_{+}$ in $S_1(x)$,
and there are $n/2$ vertices from $T_{-}$ in~$S_2(x)$.
The vertex~$z$ belongs to~$N_i$ for some~$i$. We split the case into two subcases.

(c1) If $z$ is adjacent to some vertex from~$T_{-}$, then in the neighbours
of~$x$ the function~$f$ possesses the following
values:
\begin{itemize}
    \item $f(z)=(+1,-1,0)$;
    \item the remaining $q-2$ vertices in~$N_i$ are not adjacent to any vertex from~$T_{-}$ and hence $f$ equals $(0,+1,-1)$ in them;
    \item by Lemma~\ref{l:weight1},
exactly $n$ neighbours of~$x$ are adjacent to a vertex from $T_{-}$; one of them is~$z$,
and in the remaining $n-1$ neighbours $f$ has value $(0,-1,+1)$;
    \item for all other neighbours of~$x$,
    $f$ equals $(0,0,0)$ because they are not adjacent to vertices 
    from $T_{+}\cap T_{-}$.
\end{itemize}
We summarise: $AF_x=(+1,-1,0)+(q-2) \cdot (0,+1,-1)+(n-1) \cdot (0,-1,+1)=(+1,-n+q-2, n-q+1)=(0,+1,-1)\cdot S$.

(c2) If $z$ is not adjacent to any vertex from $T_{-}$, 
then 
\begin{itemize}
    \item $f(z)=(+1,0,-1)$;
    \item by Lemma~\ref{l:weight1},
    each $N_j$, $j=1,...,n$, has exactly
    one vertex adjacent to a vertex~$v_j$ from $T_{-}$; only one of them, $v_i$,
    is adjacent to a vertex from~$T_{+}$
    and have $f(v_i)=(0,0,0)$;
    for the other $n-1$ values of~$j$, we have
    $f(v_j)=(0,-1,+1)$;
    \item the remaining $q-3$ vertices in~$N_i$,
    different from $z$ and $v_j$, the value of~$f$ is $(0,+1,-1)$;
    \item for all other neighbours of~$x$,
    not from~$N_i$ and different from~$z_j$,
    $j=1,...,n$, $f$ equals $(0,0,0)$.
\end{itemize}
In summary, 
$AF_x=            (+1, 0,-1)
     +(n-1) \cdot ( 0,-1,+1)
     +(q-3) \cdot ( 0,+1,-1)
 =(+1, -n+q-2, n-q+1)=(0,+1,-1)\cdot S$.

(d) Let $f(x)=(0,0,0)$. By definition, $x \not\in (T_{+} \cup T_{-})$. We split the case into three subcases.

(d1) For some $j$, there are $z \in T_{+}$ and $y \in T_{-}$ in $N_j$.
We have $f(z)=(+1,-1,0)$, $f(y)=(-1,+1,0)$, and $f(u)=(0,0,0)$ for each other neighbour~$u$ of~$x$. Hence, $AF_x=(0,0,0)$.

(d2) For some different $i$ and~$j$, there are $z \in T_{+}$ in~$N_i$ and $y \in T_{-}$ in~$N_j$.
We have $f(z)=(+1,-1,0)$, $f(y)=(-1,+1,0)$, 
$f(v)=(0,+1,-1)$ for each other vertex~$v$  in~$N_i$, $f(w)=(0,-1,+1)$ for each other vertex~$w$ in~$N_j$, $f(u)=(0,0,0)$ for the remaining neighbours~$u$ of~$x$. 
Again, the sum is~$(0,0,0)$.

(d3) There are no vertices from $(T_{+} \cup T_{-})$ in the neighbourhood of~$x$.
By Definition~\ref{d:weighball}, there are $l$ vertices from $T_{+}$ and $l$ vertices from $T_{-}$ in $S_2(x)$ for some~$l$.
Denote by $k$ the number of vertices adjacent
both to a vertex 
from~$T_{+}$ and to a vertex 
from~$T_{-}$. In all these vertices,
the value of~$f$ is $(0,0,0)$.
In the remaining $l-k$ neighbours of~$x$ adjacent to a vertex 
from~$T_{+}$ the value is $(0,+1,-1)$;
similarly, in the remaining $l-k$ neighbours of~$x$ adjacent to a vertex 
from~$T_{-}$ the value is $(0,-1,+1)$.
The other neighbours of~$x$ have $(0,0,0)$, and the sum of the values over the neighbourhood of~$x$ is $(0,0,0)$, which also equals~$(0,0,0)\cdot S$.

$\ref{d:matrixeq} \Longrightarrow \ref{d:weighball}$.
Let $(T_{+},T_{-})$ be an extended $1$-perfect bitrade in $H(n,q)$ according to Definition~\ref{d:matrixeq}. 
Consider an arbitrary vertex~$x$. 
Let us prove that $w_{+}(x)=w_{-}(x)$.
Recall that by definition $T_{+}$ and $T_{-}$ are codes with code distance~$4$, so 
we can use Lemma~\ref{l:weight1}. 
Consider four cases depending on the value of~$f(x)$. 

(e) Let $f(x)=(+1,-1,0)$. 
So, $x \in T_{+}$ and $w_{+}(x)=1$.
From the matrix equation, we see that the value in the first position of~$AF_x$ equals~$-1$. 
Therefore, there is a vertex from~$T_{-}$ in the neighbourhood of~$x$ and hence $w_{-}(x)=1$.

(f) Let $f(x)=(0,+1,-1)$. 
So, $x \not\in (T_{+} \cup T_{-})$ and
the value in the first position of $AF_x$ equals~$+1$.
It follows that the neighbourhood of~$x$ contains a vertex~$y$ 
from~$T_{+}$ and
does not contain any vertex
from~$T_{-}$; 
so, we have $w_{+}(x)=1$. 
It remains to show that $w_{-}(x)=1$ as well.
Let $y \in N_i$, and let $z_1,\ldots,z_k$ be the vertices 
from~$T_{-}$ in~$S_2(x)$. We consider three subcases, (f1)--(f3), but before, 
we make one useful observation.
\begin{itemize}
    \item[(*)] {\it For a vertex $v$, the second element of~$f(v)$ equals the number of neighbours of~$v$ in~$T_{+}$
    minus the number of neighbours of~$v$ in~$T_{-}$.}
    Indeed, the $v$th row of the matrix equation $AF=FS$ says that $AF_v=f(v)S$, and from the first column $(0,1,0)^{\mathrm{T}}$ of~$S$
    we see that 
    the second element of~$f(v)$
    equals  
    the first element of~$AF_v$. On the other hand, $AF_v$ is the sum of $f(z)$ over the neighbours~$z$ of~$v$.
\end{itemize}

(f1) Let $f(y)=(+1,-1,0)$.
From $AF_y=f(y)S$ (the $y$th row of the matrix equation $AF=FS$), we see that
$AF_y$ has~$-1$ in the first position; therefore, $y$ is adjacent to~$z_t$ for some~$t$.
In this case, $x$ has $q-2$ neighbours with value $(0,+1,-1)$ of~$f$
(namely, the elements $v$ of $N_i\backslash\{y\}$; by (*),
the second element of~$f(v)$ is~$+1$),
$2k-1$ vertices with value $(0,-1,+1)$ (which are adjacent to some~$z_l$; again, we use~(*)), and the remaining vertices have value $(0,0,0)$ (for the same reason).

Summing  
$f(y)=(+1,-1,0)$,
$(q-2)\cdot (0,+1,-1)$,
and $(2k-1)\cdot (0,-1,+1)$,
we get
the value $(1,-1+(q-2k-1),-(q-2k-1))$ of~$AF_x$;
on the other hand, it equals
$f(x)\cdot S = (0,+1,-1) \cdot S = (1,q-2-n,(n-1)(q-1)-n(q-2)) = (1,-1+(q-n-1),-(q-n-1)$. We deduce $k=n/2$.

(f2) Assume $f(y)=(+1,0,-1)$
and there is another vertex~$y'$
in~$N_i$ with a neighbour in~$T_{-}$.
Analogously to (f1), 
$x$ has $q-3$ neighbours with value $(0,+1,-1)$ (namely, the elements of $N_i\backslash\{y,y'\}$)
and $2k-1$ vertices with value $(0,-1,+1)$;
the remaining vertices, including~$y'$, have $(0,0,0)$.

Similarly to (f1), we derive
$AF_x=(1,-1+(q-2k-1),-(q-2k-1))$
and deduce $k=n/2$.

(f3) Assume $f(y)=(+1,0,-1)$
and there are no vertices 
in~$N_i$ with a neighbour in~$T_{-}$.
Analogously to (f2), 
$x$ has  
$q-2$ neighbours with value $(0,+1,-1)$ 
(the elements of $N_i\backslash\{y\}$)
and $2k$ vertices with value $(0,-1,+1)$.
Again, $k=n/2$. 

In all three subcases, we have $k=n/2$ and hence $w_{-}(x)=1$.

(g) Let $f(x)=(+1,0,-1)$. 
So, $w_{+}(x)=1$ and $x$ has no neighbours in~$T_{-}$.
By (*), in the neighbours of~$x$, $f$ has values $(0,0,0)$ (in the vertices that has a neighbour in~$T_-$) or $(0,+1,-1)$ (in the other vertices). 
The sum $AF_x$ of these values is equal to 
$f(x)\cdot S = (0,n(q-2),-n(q-2))$. 
Hence the neighbourhood of~$x$ contains $n$ vertices that has a neighbour in~$T_{-}$.
It follows that $|S_2(x)\cap T_{-}|=n/2$ and hence
$w_{-}(x)=1$.  

(h) Let $f(x)=(0,0,0)$. 
If there is a vertex from~$T_{+}$ and a vertex from~$T_{-}$ in the neighbourhood of~$x$, then $w_{-}(x)=w_{+}(x)=1$.
Otherwise, the number of vertices with value $(0,+1,-1)$ is equal to
the number of vertices with value $(0,-1,+1)$. 
All remaining vertices have value $(0,0,0)$.
Therefore, $S_2(x)$ has the same number of vertices from~$T_{+}$ and~$T_{-}$, i.e., $w_{+}(x)=w_{-}(x)$.
\end{proof}

The following fact is well known and straightforward.

\begin{lemma}\label{l:perfbitreigen}
A pair $(T_{+},T_{-})$ of disjoint subsets of $\VVH{n}{q}$ is a $1$-perfect bitrade in $H(n,q)$ if and only if
$T_{+}$ and $T_{-}$ are codes with code distance~$3$ and
the characteristic function~$f$ of~$(T_{+},T_{-})$ belongs to $U_l(n,q)$, where
$\lambda_l(n,q)=-1$ and, respectively, $l=\frac{n(q-1)+1}{q}$. 
\end{lemma}

\begin{proposition}\label{p:24}
 Definitions~\ref{d:proection} and~\ref{d:eigensubspace} are equivalent.
\end{proposition}
\begin{proof}
$\ref{d:eigensubspace} \Longrightarrow \ref{d:proection}$
Let $(T_{+},T_{-})$ be an extended $1$-perfect bitrade in $H(n,q)$ according to Definition~\ref{d:eigensubspace}, 
and let $f$ be its characteristic function.
By definition, $f \in U_n(n,q) \oplus U_l(n,q)$, where $\lambda_l(n,q)=q-2$, i.e., $l=\frac{n(q-1)-q+2}{q}$.
Therefore, $f=g+h$, where $g \in U_n(n,q)$, $h \in U_l(n,q)$.
Let $i \in \{1,\ldots,n\}$ be an arbitrary position. 
By Lemma~\ref{l:fiproekt}(1,2) we have $\prj{f}{i} = \prj{g}{i} + \prj{h}{i} = \prj{h}{i} \in U_{l}(n-1,q)$.
By Lemma~\ref{l:perfbitreigen} $\prj{f}{i}$ is the characteristic function of some $1$-perfect bitrade in $H(n-1,q)$.
On the other hand, the pair $(T'_{+},T'_{-})$ whose characteristic function is $\prj{f}{i}$
is the projection of $(T_{+}, T_{-})$.
Therefore, $(T_{+}, T_{-})$ is an extended $1$-perfect bitrade by Definition~\ref{d:proection}.

$\ref{d:proection} \Longrightarrow \ref{d:eigensubspace}$
Let $(T_{+},T_{-})$ be an extended $1$-perfect bitrade in $H(n,q)$ according to Definition~\ref{d:proection}, and 
let $f$ be its characteristic function.
Let $f \in U_{j_1}(n,q) \oplus \ldots \oplus U_{j_k}(n,q)$.
So, $f=g_1+\ldots+g_k$, where $g_t \in U_{j_t}(n,q)$.
Denote $l=\frac{(n-1)(q-1)+1}{q}$ (i.e., $\lambda_l(n-1,q)=-1$ in $H(n-1,q)$ for this $l$).
Suppose that there is $m$ in $\{j_1,\ldots,j_k\} \backslash \{l,n\}$.
By Lemma~\ref{l:fiproekt}(3), there is~$r$ such that 
$\prj{g_m}{r} \not\equiv 0$.
Hence, $\prj{f}{r} = \prj{g_1}{r} + \ldots + \prj{g_k}{r}=h_1+\ldots+h_k$, 
where $h_t \in U_{j_t}(n-1,q)$ and $h_m \not\equiv 0$. Therefore, $\prj{f}{r} \not\in U_{l}(n-1,q)$, 
and we have a contradiction with the fact that $\prj{f}{r}$ is the characteristic function of some $1$-perfect bitrade.
\end{proof}

Propositions~\ref{p:235}--\ref{p:24} 
prove Theorem~\ref{th:equiv}.

A pair $(T_{+},T_{-})$ of disjoint subsets of $\VVH{n}{q}$ is called an \emph{extended $1$-perfect bitrade} in $H(n,q)$ 
if it satisfies one of equivalent Definitions~\ref{d:first}--\ref{d:last}. 

\section{Existence of extended 1-perfect bitrades}\label{s:necessary}

To construct extended $1$-perfect bitrades, we adopt the approach from~\cite{MogSol:perfbit},
which solves the problem of existence of $1$-perfect bitrades.
In~\cite{MogSol:perfbit},
objects called spherical bitrades are considered. 
To adopt the technique to extended $1$-perfect bitrades,
we need spherical bitrades with an additional condition, extendability.

Recall that a code~$C$ with distance~$d$ in $H(n,q)$ 
is called an \emph{MDS code} 
if $|C|=q^{n-d+1}$. 
MDS codes have the following useful property:
if $C$ is an MDS code, then a punctured code (obtained by deleting some coordinate) and a shortened code (the part of the punctured code with a fixed symbol in the deleted coordinate) are MDS too.

\begin{definition}\label{d:sph} 
A pair $(T_{+},T_{-})$ of disjoint subsets of $\VVH{n}{q}$
is called an \emph{extendable spherical bitrade} in $H(n,q)$ 
if 
\begin{enumerate}
    \item the characteristic function of $(T_{+},T_{-})$ is an eigenfunction with eigenvalue $0$;
    \item $T_{+}$ and $T_{-}$ are codes with code distance~$3$;
    \item $T_{+}$ and $T_{-}$ can be represented as $T_{+}=T^0_{+} \cup \ldots \cup T^{q-1}_{+}$ and $T_{-}=T^0_{-} \cup \ldots \cup T^{q-1}_{-}$, where the union is disjoint, $T^i_{+}$ and $T^i_{-}$ are codes with code distance~$4$, and $T^i_{+} \cup T^i_{-}$ is a code with code distance~$3$,
    $i =0,1,\ldots,q-1$.
    \end{enumerate}
\end{definition}

\begin{lemma}
\label{l:sperbitind}
    If $(T_{+},T_{-})$ is an extendable spherical bitrade, then $T_{+} \cup T_{-}$ is an independent set (i.e., a distance-$2$ code).
\end{lemma} 
\begin{proof}
    Condition 1 of Definition~\ref{d:sph} implies
    that $|S_1\cap T_+|=|S_1\cap T_-|$ for any vertex~$x$. If $x \in T_+$ (similarly, for $x \in T_-$), then from condition 2 we have $|S_1\cap T_+|=0$ and hence $|S_1\cap T_-|=0$.
\end{proof}

\begin{lemma}\label{l:qq}
If $q=2^m$, $m\ge 2$, then there is an extendable spherical bitrade in $H(q,q)$.
\end{lemma}
\begin{proof}
Let $C$ be an extended $1$-perfect code in $H(q+2,q)$. 
It is also an MDS code with code distance~$4$. 
Consider the following sets of vertices in $\VVH{q}{q}$:
$$C_{i,j}=\{(x_1,\ldots,x_q) \in \VVH{q}{q}:\  (x_1,\ldots,x_q,i,j) \in C\},$$ $$C_{i;}=C_{i,0} \cup \ldots \cup C_{i,q-1} \quad\mbox{and}\quad C_{;j}=C_{0,j} \cup \ldots \cup C_{q-1,j}.$$ 
From the properties of MDS codes, we have that $C_{i,j}$ 
is an MDS code with code distance~$4$,
$C_{i;}$ and $C_{;j}$ are MDS codes with code distance~$3$,  
and $C'=C_{0,0} \cup \ldots \cup C_{q-1,q-1}$ 
is an MDS code with code distance~$2$. 
Define $T_{+}=C_{;0}$ and $T_{-}=C_{;1}$. 
Also, define $T^i_{+}=C_{i,0}$ and $T^i_{-}=C_{i,1}$. 
Let us prove that $(T_{+}, T_{-})$ 
is an extendable spherical bitrade. Indeed, 
$T^i_{+} \cup T^i_{-}=C_{i,0} \cup C_{i,1}$ 
is a code with code distance~$3$; the remaining requirements of
p.2 and p.3 hold from the definition. It remains to prove p.1. Let $f$ be the characteristic function of $(T_{+},T_{-})$. If $x\in C'$, then $f(y)=0$ for each neighbour~$y$ of~$x$; so, the eigenfunction condition   holds in~$x$. Suppose $x\not\in C'$. Since $C'$ is an MDS code with code distance~$2$, $x$ has $q$ neighbours in~$C'$.
On the other hand, $C'$ is divided into $q$ codes with code distance~$3$, $C_{;0}, \ldots, C_{;q-1}$; so, $x$ has exactly one neighbour in~$C_{;0}$ and one in~$C_{;1}$. Therefore, the sum of the values of~$f$ over the neighbourhood of~$x$ equals~$0$. By the definition, $f$ is an eigenfunction with eigenvalue~$0$. 
\end{proof}

For two functions $g:\VVH{n_1}{q} \to \RR$ and $h:\VVH{n_2}{q} \to \RR$,
their \emph{tensor product} $f=g \otimes h:\VVH{n_1+n_2}{q} \to \RR$ is defined
as follows:
$f(x,y)=g(x)\cdot h(y)$, where $x \in \VVH{n_1}{q}$ and $y \in \VVH{n_2}{q}$. 

\begin{lemma}[\cite{ValVor:Hamm}]\label{l:oplus}
If $g$ is an eigenfunction of $H(n_1,q)$ with eigenvalue~$\lambda$ and $h$ is an eigenfunction of $H(n_2,q)$ with eigenvalue~$\mu$, then $g \otimes h$ is an eigenfunction of $H(n_1+n_2,q)$ with eigenvalue $\lambda+\mu$.
\end{lemma}

\begin{lemma}\label{l:sph}
If $q=2^m\ge 4$, then there is an extendable spherical bitrade in $H(lq,q)$ for any positive integer~$l$.
\end{lemma}
\begin{proof}
Let us prove the statement by induction on~$l$. For $l=1$, it is true by Lemma~\ref{l:qq}. 
Suppose that there is an extendable spherical bitrade $(A_{+}, A_{-})$ in $H(lq,q)$. Denote its characteristic function by~$f_A$.  Let us construct an extendable spherical bitrade in $H(q(l+1),q)$. Consider some extendable spherical bitrade in $H(q,q)$ and denote it by $(B_{+},B_{-})$. Let $f_B$ be its characteristic function. Consider the vertices of $H(q(l+1),q)$ as $(x,y)$, where $x \in \VVH{lq}{q}$ and $y \in \VVH{q}{q}$. 

The tensor product $f_A \otimes f_B$
is the characteristic function
of some pair $(T_{+}, T_{-})$
of disjoint sets. 
Let us prove that $(T_{+}, T_{-})$
 is an extendable spherical bitrade in $H(q(l+1),q)$.

1) As $f_A \otimes f_B$ is the tensor product of eigenfunctions with eigenvalue~$0$,
by Lemma~\ref{l:oplus} it is itself an eigenfunction with eigenvalue~$0$.

2) The code distance of $T_{+}$ is~$3$. Indeed, $T_{+}=A_{+} \times B_{+} \cup A_{-}\times B_{-}$. Let $(x,y)$ and $(x',y')$ be two different vertices from~$T_{+}$. If $x=x'$, then either $y, y' \in B_{+}$ or $y, y' \in B_{-}$; therefore, the distance between $(x,y)$ and $(x',y')$ is not less than~$3$. The case $y=y'$ is similar. If $x \ne x'$ and $y \ne y'$, then by Lemma~\ref{l:sperbitind} we have $d(x,x')\ge 2$ and $d(y,y')\ge 2$,
and hence $d((x,y),(x',y'))\ge 4$. Analogously, $T_{-}$ has code distance~$3$.

3) By the induction assumption, $$A_{+}=A^0_{+} \cup \ldots \cup A^{q-1}_{+}, A_{-}=A^0_{-} \cup \ldots \cup A^{q-1}_{-},$$ $$B_{+}=B^0_{+} \cup \ldots \cup B^{q-1}_{+}, B_{-}=B^0_{-} \cup \ldots \cup B^{q-1}_{-},$$ 
where
$A^i_{+}$, $A^i_{-}$, $B^i_{+}$, $B^i_{-}$ are codes with code distance~$4$ and $A^i_{+} \cup A^i_{-}$ and $B^i_{+} \cup B^i_{-}$ are codes with code distance~$3$ for every $i \in \{0,1, \ldots, q-1\}$.
Note that $T_{+} \cup T_{-}$ consists of all sets of form $A^i_{\delta} \times B^j_{\sigma}$, where $\delta, \sigma \in \{+,-\}$. 
Define a partition $T_{+}=T^0_{+} \cup \ldots \cup T^{q-1}_{+}$ and $T_{-}=T^0_{+} \cup \ldots \cup T^{q-1}_{+}$ in the following way: $A^i_{\delta}\times B^j_{\sigma}$, where $\delta, \sigma \in \{+,-\}$, belongs to $T^{i+j \mod q}_\theta$, where $\theta$ is $+$ if $\delta=\sigma$, and $\theta$ is $-$ if $\delta \ne \sigma$. 
Let $u=(x,y)$ and $v=(x',y')$ be vertices from some $T^s_{\gamma}$.
Suppose $u$ belongs to some $A^i_{\sigma} \times B^k_{\theta}$ and 
$v$ belongs to some $A^j_{\delta} \times B^l_{\eta}$. 
If  $A^i_{\sigma}$ and $A^j_{\delta}$  coincide, then  $B^j_{\theta}$ and $B^l_{\eta}$  coincide by the definition.  
In this case, $d(u,v) \ge 4$ because $u$ and~$v$ belong to the direct product of codes with code distance~$4$. 
If $B^k_{\theta}$ and $B^l_{\eta}$  coincide, then analogously $d(u,v) \ge 4$. Let $i \ne j$, and hence $k \ne l$. 
In this case, the distance between~$x$ and~$x'$ in $H(lq,q)$ is not less than~$2$, and 
the distance between~$y$ and~$y'$ in $H(q,q)$ is not less than~$2$. 
Hence, $d(u,v) \ge 4$. Let $u=(x,y)$ belong to $T^s_{+}$ and $v=(x',y')$ to $T^s_{-}$ for some~$s$. 
Let $x \in A^i_{\sigma}$, $y \in B^k_{\theta}$, $x' \in A^j_{\delta}$, $y \in B^l_{\eta}$. 
If $i \ne j$, then $k \ne l$; hence the sets~$A^i_{\sigma}$ and~$A^j_{\delta}$ do not coincide, and the sets $B^k_{\theta}$ and $B^l_{\eta}$ are also different. Therefore, the code distance between~$u$ and~$v$ is at least~$4$. 
Let $i=j$ and $k=l$. Suppose $\sigma=\delta$ (the case $\theta=\eta$ is similar).
Then the distance between~$u$ and~$v$ not less than the distance between~$y$ and~$y'$ in $H(q,q)$ and cannot be less than~$3$.
Therefore, 
p.3 also holds and $(T_{+}, T_{-})$ is an extendable spherical bitrade.
\end{proof}

\begin{definition}\label{d:pb}
 A pair $(T_{+},T_{-})$ of disjoint subsets of $\VVH{n}{q}$
is called an \emph{extendable $1$-perfect bitrade} in $H(n,q)$ 
if 
\begin{enumerate}
    \item the characteristic function of $(T_{+},T_{-})$ is an eigenfunction with eigenvalue $-1$;
    \item $T_{+}$ and $T_{-}$ are codes with code distance $3$;
    \item $T_{+}$ and $T_{-}$ can be represented as $T_{+}=T^0_{+} \cup \ldots \cup T^{q-1}_{+}$ and $T_{-}=T^0_{-} \cup \ldots \cup T^{q-1}_{-}$, where $T^i_{+}$ and $T^i_{-}$ are codes with code distance~$4$ for every 
    $ i \in \{0,1,\ldots,q-1\}$.
    \end{enumerate}   
\end{definition}

\begin{remark}
As we will see in the proof of Theorem~\ref{t:exofbitrade} below,
an extendable $1$-perfect bitrade
can be turned to an extended $1$-perfect bitrade from equivalent 
Definitions~\ref{d:proection}--\ref{d:eigensubspace}
by appending the symbol~$i$ 
to every word from $T^i_{+} \cup T^i_{i}$.
The inverse is also true:
if we have an extended $1$-perfect bitrade,
then puncturing (removing the last symbol
from all its words) results in an  extended $1$-perfect bitrade.
\end{remark}

\begin{lemma}\label{l:pb}
    If $q=2^m\ge 4$, then for every positive integer~$l$ there is an extendable $1$-perfect bitrade in $H(lq+1,q)$.
\end{lemma}
\begin{proof}
    By Lemma~\ref{l:sph}, there is an extendable spherical bitrade $(A_{+},A_{-})$.
    Let $f_0$ be its characteristic function, and let
    $A_{+}=A^0_{+} \cup \ldots \cup A^{q-1}_{+}$, $A_{-}=A^0_{-} \cup \ldots \cup A^{q-1}_{-}$,
    according to Definition~\ref{d:sph}. Consider the function~$f_1$ on $H(1,q)$ defined as follows: $f_1(0)=1$, $f_1(1)=-1$, and $f_1(x)=0$ if $x \ne 0,1$. Clearly, $f_1$ is an eigenfunction with eigenvalue~$-1$. Consider the function $f=f_0 \otimes f_1$. Define the pair $(T_{+}, T_{-})$ such that $f$ is its characteristic function. Let us prove that $(T_{+},T_{-})$ is an extendable $1$-perfect bitrade. Consider the vertices of $H(lq+1,q)$ as $(x,y)$, where $x \in \VVH{lq}{q}$ and $y \in \VVH{}{q}$.
    
    1) At first, $f$ is an eigenfunction with eigenvalue~$-1$ as the tensor product of the eigenfunctions $f_0$ and $f_1$. 
    
    2) At second, $T_{+}$ and $T_{-}$ are codes with code distance~$3$. Indeed, let $u=(x,y)$ and $v=(x',y')$ be vertices from~$T_{+}$. The distance between~$x$ and~$x'$ in $H(lq,q)$ is not less than~$2$ (if $x=x'$, then $y=y'$) and, moreover, not less~$3$ if $y=y'$ (since in this case~$x$ and~$x'$ either both belong to~$A_{+}$ or both belong to $A_{-}$). Therefore, $T_{+}$ (analogously $T_{-}$) is a code with distance~$3$.
    
    3) Define $T^i_{+}=A^i_{+} \times \{0\} \cup A^i_{-} \times \{1\}$ and 
    $T^i_{-}=A^i_{-} \times \{0\} \cup A^i_{+} \times \{1\}$. Let $u=(x,y)$ and $v=(x',y')$ be vertices from~$T^i_{+}$ for some~$i$. The distance between~$x$ and~$x'$ is not less than~$3$ and, moreover, not less than~$4$ if $y=y'$ (since in this case either $x$ and $x'$ both belong to $A^i_{+}$ or they both belong to $A^i_{-}$). Therefore, $T^i_{+}$ (analogously $T^i_{-}$) is a code with distance~$4$ for each $i \in \{0,1, \ldots,q-1\}$.       
\end{proof}

\begin{theorem}\label{t:exofbitrade}
    If $q=2^m$ then for any positive integer~$l$ there is an extended $1$-perfect bitrade in $H(lq+2,q)$.
\end{theorem}
\begin{proof}
    By Lemma~\ref{l:pb}, there is an extendable $1$-perfect bitrade in $H(lq+1,q)$, say $(A_{+},A_{-})$.  
    Let $A_{+}=A^0_{+} \cup \ldots \cup A^{q-1}_{+}$, $A_{-}=A^0_{-} \cup \ldots \cup A^{q-1}_{-}$ be its representation according to Definition~\ref{d:pb}. Let us define $(T_{+},T_{-})$ in the following way:
    $T_{+}=T^0_{+} \cup \ldots \cup T^{q-1}_{+}$ and  $T_{-}=T^0_{-} \cup \ldots \cup T^{q-1}_{-}$, where $T^i_{+}=A^i_{+}\times \{i\}$, $T^i_{-}=A^i_{-} \times \{i\}$. Let us prove that $(T_{+},T_{-})$ is an extended $1$-perfect bitrade in $H(lq+2,q)$. Indeed, let $u=(x, y)$ and $v=(x',y')$ be vertices from~$T_{+}$, where $x \in \VVH{lq+1}{q}$ and $y \in \VVH{}{q}$. The distance between~$x$ and~$x'$ in $H(lq+1,q)$ is not less than~$3$ and not less than $4$ if $y=y'$ (since in this case, $x$ and $x'$ belong to $A^i_{+}$). Therefore, $T_{+}$ (similarly,~$T_{-}$) is a code with distance~$4$. Let $f$ be the characteristic function of $(T_{+},T_{-})$. By the construction, the $(lq+2)$-projection of $(T_{+},T_{-})$ is $(A_{+},A_{-})$, with characteristic function  $\prj{f}{lq+2}$.  Since the characteristic function of $(A_{+},A_{-})$ is an eigenfunction with eigenvalue~$-1$,
    the function~$f$ belongs to the direct sum of the eigenspaces corresponding to the eigenvalues~$q-2$ and~$-n$. Therefore, $(T_{+},T_{-})$ is an extended $1$-perfect bitrade by Definition~\ref{d:eigensubspace}.
\end{proof}

\begin{theorem}\label{th:necessary}
If there is an extended $1$-perfect bitrade $(T_{+},T_{-})$ in $H(n,q)$, then 
\begin{enumerate}
    \item $n$ is even,
    \item $n=lq+2$ for some $l \in \{0,1,\ldots\}$.
\end{enumerate}
\end{theorem}
\begin{proof}
1) Let $x \in T_{+}$. 
Then there is a vertex $y \in S_1(x)$ such that $y$ is not adjacent to any vertex from $T_{-}$.
Since $w_{+}(y)=1$, there are $n/2$ vertices from $T_{-}$ in $S_2(y)$.
Hence, $n$ is even.

2) By Definition~\ref{d:eigensubspace}, $n(q-1)-iq$ is $q-2$ for some~$i$. 
Hence, $n\equiv2 \mod q$.
\end{proof}

\begin{corollary}
Extended $1$-perfect bitrades in $H(n,q)$, $q=2^m$, exist if and only if $n=lq+2$ for some positive integer~$l$.
\end{corollary}

For a vertex $v$ of ${H(n,q)}$,
denote by $W^{+}_{v,i}$ and $W^{-}_{v,i}$ 
the number of vertices 
at distance~$i$ from~$v$
that belong to~$T_{+}$ 
and~$T_{-}$, respectively.

\begin{theorem}\label{t:weightdist}
Let $(T_{+},T_{-})$ be an extended $1$-perfect bitrade in $H(n,q)$. 
Then for any vertex~$v$ and $i \in \{1,\ldots,n\}$ the value $W^{+}_{v,i}-W^{-}_{v,i}$ depends only 
on~$i$ and values $W^+_{v,0}$, $W^{-}_{v,0}$, $W^{+}_{v,1}$ and $W^{-}_{v,1}$.
\end{theorem}
\begin{proof}
By Definition~\ref{d:matrixeq} there is a matrix $F$
whose first column $F_{0}$ is the
characteristic function of 
$(T_{+},T_{-})$ 
such that $AF=FS$.
This matrix equation gives $A^kF=FS^k$ for any $k$; hence, $P(A)F=FP(S)$ for any polynomial $P$.
It is known (see, e.g.,~\cite[Proposition~2.7.1]{Brouwer}) that $A_l=P_l(A)$ for some polynomial $P_l$, where $A_l$ is the distance-$l$ adjacency matrix
(i.e., for any vertices $u$ and $t$ the $(u,t)$-th element of $A_l$ equals $1$ if $d(u,t)=l$, and equals $0$ otherwise).

By the definition,
$W^{+}_{v,i}-W^{-}_{v,i}$
is the sum of the
characteristic function of
$(T_{+},T_{-})$
over the vertices $u$ at 
distance $i$ from $v$.
Trivially,
it is the $v$th element
of $A_i F_{0}$,
where $F_0$ 
is the first column of $F$ 
and the characteristic vector
of $(T_{+},T_{-})$.
Since $A_iF=P_i(A)F=FP_i(S)$, the value $W^{+}_{v,i}-W^{-}_{v,i}$ depends only on $F^{(v)}$, 
where $F^{(v)}$ is $v$-row of $F$.
It remains to observe
that $F^{(v)} = (a_0,a_1,a_2)$
where

$a_0 = W^{+}_{v,0}-W^{-}_{v,0}$, by (i) in Definition~\ref{d:matrixeq},

$a_1 = W^{+}_{v,1}-W^{-}_{v,1}$,
from (iii) in
Definition~\ref{d:matrixeq}
(since $a_1$ equals the $(v,1)$-element of $FS$ and $W^{+}_{v,1}-W^{-}_{v,1}$ equals the $(v,1)$-element of $AF$),

$a_2 = -a_0-a_1$,
by (ii) in
Definition~\ref{d:matrixeq}.
\end{proof}

The polynomial from Theorem~\ref{t:weightdist} is called \emph{Krawtchouk polynomial}. 
The statement of the theorem gives a possible way to reject some putative parameters of extended $1$-perfect bitrades 
by calculating the values $W^{+}_{v,i}-W^{-}_{v,i}$. 
If at least one of them is noninteger, then there are no extended perfect bitrades for given~$n$ and~$q$.
This method is well known for completely regular codes and equitable partitions, see for example~\cite{Kro:struct}.

\section{Conclusion}\label{s:conc}
Bitrades (trades) can be defined
for different classes
of combinatorial objects,
for example, combinatorial designs
\cite{HedKho:trades} and their subspace
analogues~\cite{KMP:16:trades},
latin squares~\cite{Cav:rev} and 
latin hypercubes~\cite{Potapov:2013:trade}. 
Trades help to construct objects of the corresponding class and the study of trades can potentially help to prove the nonexistence of objects with certain parameters.
The spectrum of volumes of trades is
related to possible differences between
two objects from the corresponding class 
on the one hand and to possible sizes
of the support for functions with spectrum constraint \cite{SotVal:survey2021} on the other hand.
We have shown how to define
bitrades for diameter perfect codes,
for uniformly packed codes,
and for completely regular codes,
and proved that these definitions
are equivalent in the case
of extended $1$-perfect codes
(defined as a diameter perfect, 
uniformly packed, 
or completely regular code,
or in two other manners).
We consider these results
as a contribution to the general theory of trades,
while the results of Section~\ref{s:necessary}
are more special, for extended $1$-perfect codes.

\subsection*{Acknowledgment}
The results in Sections~\ref{s:def}, \ref{s:equivalence} and Theorem~\ref{th:necessary} in Section~\ref{s:necessary} were obtained at the expense of the Russian Science Foundation (Grant 18-11-00136), \url{https://rscf.ru/en/project/18-11-00136/}; the other results in Section~\ref{s:necessary} were obtained within the framework
of the state contract of the Sobolev Institute of Mathematics (FWNF-2022-0017).

\subsection*{Declaration of competing interest}
The authors declare that they have no known competing financial interests or personal
relationships that could have appeared to influence the work reported in this paper.

\subsection*{Data availability}
No data was used or generated during the research described in the article.



\providecommand\href[2]{#2} \providecommand\url[1]{\href{#1}{#1}}
  \def\DOI#1{{\href{https://doi.org/#1}{https://doi.org/#1}}}

\end{document}